\documentclass[10pt,twocolumn,twoside]{IEEEtran} 	
\usepackage{amsmath,amssymb,bm,bbm,mathrsfs,amscd,color,hyperref}
\usepackage{amsthm,subfigure}
\usepackage{dsfont, color}
\usepackage{graphicx}
\usepackage{epsfig}
\def\qed{\hfill \vrule height 7pt width 7pt depth 0pt\medskip}
\def\beq{\begin{equation}}
\def\eeq{\end{equation}}
\def\proof{\noindent{\bf Proof}\ \ }
\newtheorem{theorem}{Theorem}

\newtheorem{proposition}[theorem]{Proposition}
\newtheorem{lemma}[theorem]{Lemma}

\newtheorem{example}{Example}
\theoremstyle{remark}

\newcommand{\ds}{\displaystyle}

\newcommand{\ba}{\begin{array}}
\newcommand{\ea}{\end{array}}

\renewcommand{\l}{\left}\renewcommand{\r}{\right}

\newcommand{\be}{\begin{equation}}
\newcommand{\ee}{\end{equation}}
\newcommand{\eps}{\varepsilon}

\newcommand{\mc}{\mathcal}

\newcommand{\ov}{\overline}

\newcommand{\Z}{\mathbb{Z}}

\newcommand{\1}{\mathbbm{1}}
\newcommand{\E}{\mathbb{E}}

\newcommand{\R}{\mathbb{R}}

\renewcommand{\P}{\mathbb{P}}

\DeclareMathOperator{\supp}{supp}

\DeclareMathOperator{\tmix}{\tau_{\mathrm{mix}}}

\def\1{\mathds{1}}
\def\Z{\mathbb{Z}}

\def\E{\mathbb{E}}
\def\R{\mathbb{R}}

\def\P{\mathbb{P}}

\begin{document}


\title{Robustness of large-scale stochastic matrices to localized perturbations}



\author{Giacomo Como and Fabio Fagnani
\thanks{ G.~Como is with the Department of Automatic Control, Lund University, SE-221 00 Lund, Sweden {\tt\small giacomo.como@control.lth.se}. F.~Fagnani is with the Department of Mathematical Sciences, Politecnico di Torino, Torino, 10129 {\tt\small fabio.fagnani@polito.it}. }}%








\maketitle
\begin{abstract}
Many notions of network centrality can be formulated in terms of invariant probability vectors of suitably defined stochastic matrices encoding the network structure. Analogously, invariant probability vectors of stochastic matrices allow one to characterize the asymptotic behavior of many linear network dynamics, e.g., arising in opinion dynamics in social networks as well as in distributed averaging algorithms for estimation or control.  Hence, a central problem in network science and engineering is that of assessing the robustness of such invariant probability vectors to perturbations possibly localized on some relatively small part of the network. 
In this work, upper bounds  are derived on the total variation distance between the invariant probability vectors of two stochastic matrices differing on a subset $\mc W$ of rows. Such bounds depend on three parameters: the mixing time and the entrance time on the set $\mc W$ for the Markov chain associated to one of the matrices; and the exit probability from the set $\mc W$ for the Markov chain associated to the other matrix. These results, obtained through coupling techniques, prove particularly useful in scenarios where $\mc W$ is a small subset of the state space, even if the difference between the two matrices is not small in any norm. Several applications to large-scale network problems are discussed, including robustness of Google's PageRank algorithm, distributed averaging, consensus algorithms, and the voter model. 

\end{abstract}


\section{Introduction}
How much can the invariant probability vector
$$\pi=\pi P$$ of an irreducible row-stochastic matrix $P$ be affected by perturbations localized on a relatively small subset $\mc W$ of its state space $\mc V$? 
Such a question arises in an increasing number of applications, most notably in the emerging field of large-scale networks. 

As an example, many notions of network centrality  can be formulated in terms of invariant probability vectors of suitably defined stochastic matrices.
In particular, Google's PageRank algorithm \cite{BrinPage:98} assigns to webpages values corresponding to the entries of the invariant probability vector $\pi$ of the matrix $P$ obtained as a convex combination of the normalized adjacency matrix of the directed graph describing the hyperlink structure of the World Wide Web (WWW), and of a matrix whose entries are all equal to the inverse of the total number of webpages \cite{LangvilleMeyer:2006,Chung:2007}. 
A well-known problem in this context is rank-manipulation, i.e., the intentional addition or removal of hyperlinks from some webpages (hence, the alteration of the corresponding rows of $P$) with the goal of modifying the PageRank vector \cite{AvrachenkovLitvak:2006,KerchoveNinoveVanDooren:2008,CsajiJungersBlondel:2013}. A natural question is then, to what extent a small subset $\mc W$ of webpages can alter the PageRank vector $\pi$. Similar robustness issues have been raised for accidental variations of the WWW topology occurring, e.g., because of server failures or network congestion problems \cite{IshiiTempo:2009}.

The problem is of central interest also in the context of distributed averaging and consensus algorithms \cite{OlshevskyTsitsiklis:2011}. There, linear systems of the form $$x(t+1)=Px(t)\,,$$ or their continuous-time analogues, are studied, e.g., as algorithms for distributed optimization \cite{johnthes,distasyn}, control  \cite{ali,reza}, synchronization in sensor networks \cite{QunRus:2006}, or reputation management in ad-hoc networks \cite{LiuYang:2003}, as well as behavioral models for flocking phenomena \cite{vicsek}, 
 or opinion dynamics in social networks \cite{Degroot,demarzo,golub,Acemoglu.etal:2013}. Equilibria of such systems are consensus vectors, i.e., multiples of the all-one vector $\1$, and standard results following from Perron-Frobenius theory guarantee convergence 
 (under the assumption of irreducibility and, in the discrete-time case, acyclicity of $P$) 
 to a consensus vector 
 $$\ov x\1\,,\qquad \ov x=\pi x(0)\,.$$ 
Depending on the specific application, the natural question is to what extent the consensus value $\ov x$ is affected by perturbations of $P$ corresponding, e.g., to malfunctioning of a small fraction of the sensors, or conservative/influential minorities in social networks \cite{influence}.  

Other applications can be found in the context of interacting particle systems \cite{Liggettbook,Liggettbook2}. In particular, in the voter model on a finite graph \cite{ConnelyWelsh82,Cox89}, \cite[Ch.~14]{aldous-fill}, 
\cite[Ch.~6.9]{Durrettbook}, the probability vector of the final consensus value is determined by the invariant probability vector of the stochastic matrix associated to the simple random walk on the graph. Perturbations in this case may model the presence of inhomogeneities or `zealots'  \cite{Mobilia03,Mobilia05}, namely agents with an asymmetric behavior in the way they influence and are influenced from their neighbors. 

The above-described problems all boil down to estimating the distance between the invariant probability vector $\pi$ of an irreducible stochastic matrix $P$ and an invariant probability vector $\tilde\pi=\tilde\pi\tilde P$ of another stochastic matrix $\tilde P$, to be interpreted as a perturbed version of $P$. In some applications, $P$ may be reversible, equivalently be obtained by normalizing the rows of a symmetric nonnegative matrix $W$, and $\pi$ can be explicitly computed in terms of the row sums of $W$. However, even in these cases, the considered perturbations will typically be such that $\tilde P$ is not reversible and thus $\tilde\pi$ does not allow for a tractable explicit expression.

Remarkably, standard perturbation results based on sensitivity analysis \cite{Schweitzer:68,Seneta:1988,Seneta:1993,Meyer:1994,ChoMeyer:2000,ChoMeyer:2001,Mitrophanov:2003,Mitrophanov:2005,influence} do not provide a satisfactory answer to this problem. Indeed, they provide upper bounds of the form 
\be\label{eq:sensitivity}||\tilde \pi-\pi||_p\le\kappa_P ||\tilde P-P||_{q}\,,\ee
for some $p,q\in[1,\infty]$, where $\kappa_P$ is a condition number depending on the original stochastic matrix $P$ only.  
Such condition numbers are lower bounded by an absolute positive constant (e.g., $1/4$ for the smallest of those surveyed in \cite{ChoMeyer:2001}) and typically blow up as the state space $\mc V$ grows large. Therefore, such results do not allow one to prove that the distance $||\tilde\pi-\pi||_p$ vanishes in the limit of large network size, even if the stochastic matrices $P$ and $\tilde P$ differ only in a single row, unless  $||\tilde P-P||_q$ itself vanishes.

In this paper, we obtain upper bounds on the total variation distance $||\tilde\pi-\pi||:=\frac12||\tilde\pi-\pi||_1$ of the form 
\be\label{eq:mainresult}||\tilde \pi-\pi||\le\Psi\l(\frac{\tmix}{\tilde\gamma_{\mc W}\cdot\tau^*_{\mc W}}\r)\,,\ee
(see Theorem \ref{theorem:main}) where: 
\begin{itemize}
\item $\Psi:[0,+\infty)\to[0,1]$ is a continuous, nondecreasing function such that $\Psi(0)=0$ (see \eqref{eq:thetadef} for its definition and Figure \ref{fig:thetaplot} for its graph);
\item $\tmix$ is the mixing time of the matrix $P$, defined as
\be\label{eq:mixingdef}\tmix:=\inf\l\{t\ge1:\ \max_{u,v\in\mc V}||P^t_{u,\cdot}-P^t_{v,\cdot}||\le\frac1e\r\}\,,\ee
i.e., as the minimum $t$ such that all the rows of the $t$-th power $P^t$ are within total variation distance $1/e$ from each other;  
\item $\tau^*_{\mc W}$ is the entrance time on the set $\mc W$, 
defined as 
\be\label{eq:tauWdef}\tau^*_{\mc W}:=\min_{u\in\mc V\setminus\mc W}\tau^u_{\mc W}\,\ee
where $\tau^u_{\mc W}$, for $u\in\mc V$, are the solution of the linear system
\be\label{hittingdef}
\tau^u_{\mc W}=0\,,\ u\in\mc W\,,\qquad \tau^u_{\mc W}=1+\sum_{v\in\mc V}P_{uv}\tau^v_{\mc W}\,,\ u\in\mc V\setminus\mc W
\ee
and thus coincide with the expected hitting times on the set $\mc W$ for a Markov chain with transition probability matrix $P$;
\item $\tilde\gamma_{\mc W}$ stands for the exit probability  from $\mc W$ defined as 
\be\label{def:exitprobability}
\tilde\gamma_{\mc W}:=\sup_{t\ge0}\min_{\substack{w\in\mc W:\\\tilde\pi_w>0}}\frac1t\sum_{k=1}^t\ 
\sum_{\substack{\xi_0=w,\xi_k\in\mc V\setminus\mc W\\\xi_1,\ldots, \xi_{k-1}\in\mc W}}\ \prod_{1\le l\le k}\tilde P_{\xi_{l-1}\xi_l}\,,
\ee
where the second summation runs over all $(k+1)$-tuples $\xi$ that start with $\xi_0=w$, end with some $\xi_{k}\in\mc V\setminus\mc W$, and have all intermediate entries $\xi_l$ in $\mc W$, for $1\le l<k$. As shown in 
\eqref{eq:kappaWdef}, the argument of the minimization in \eqref{def:exitprobability} coincides with the probability that a Markov chain with transition probability matrix $\tilde P$ started at $w$ exits from $\mc W$ before time $t$, normalized by $t$.
\end{itemize}

The bound \eqref{eq:mainresult} is proved in Section \ref{sect:mainresult} and constitutes the main technical contribution of this paper. To the best of our knowledge, this result is completely original. In particular, it does not follow from the aforementioned sensitivity results. As opposed to them, whose proofs are all algebraic in nature, its proof relies on probabilistic coupling techniques, combined with an argument similar to the one developed in \cite{Acemoglu.etal:2013} in the context of `highly fluid' social networks. 
Because of the properties of the function $\Psi(\,\cdot\,)$, the bound (\ref{eq:mainresult}) implies that the total variation distance $||\tilde\pi-\pi||_{}$ vanishes provided that the quantity $\tmix/(\tilde\gamma_{\mc W}\cdot\tau^*_{\mc W})$ does. 
As we will show, this finds immediate application in the PageRank manipulation problem. 
More in general, our results prove useful in many of those large-scale network applications where classical sensitivity-based results fail to provide a satisfactory answer. 

Mixing properties of stochastic matrices have been the object of extensive recent research \cite{aldous-fill,MontenegroTetali,LevinPeresWilmer}, and several results are available allowing one to estimate the mixing time $\tmix$ of a stochastic matrix $P$, e.g., in terms of the conductance or other geometrical properties of the graph associated to $P$. It is worth pointing out that a connection between mixing properties and robustness of stochastic matrices is already unveiled by the perturbation results of \cite{Mitrophanov:2003,Mitrophanov:2005}, where (\ref{eq:sensitivity}) is proven for $p=1$, $q=\infty$, and condition number $\kappa_P$ proportional to $\tmix$. Of a similar flavor are Seneta's results \cite{Seneta:1988,Seneta:1993} estimating the condition number $\kappa_P$ in terms of ergodicity coefficients. Also the estimates proposed in \cite{influence} for symmetric $P$ can be rewritten as (\ref{eq:sensitivity}) with for $p=q=2$ and $\kappa_P$ equal to the inverse of the spectral gap of $P$. As compared to these references, the fundamental novelty of our bound (\ref{eq:mainresult}) consists in measuring the size of the perturbation in terms of $1/(\tilde\gamma_{\mc W}\cdot\tau^*_{\mc W})$ instead of the distance $||\tilde P-P||_q$, thus enabling one to obtain significant results in scenarios where $\mc W$ is small but $\tilde P-P$ is not necessarily small in any norm.

In fact, of the parameters appearing in the righthand side of \eqref{eq:mainresult},  the exit probability $\tilde\gamma_{\mc W}$ is the only one truly depending on the perturbation $\tilde P-P$, and is indeed easily estimated in typical cases when $\mc W$ is a small subset of $\mc V$. On the other hand, the entrance time $\tau^*_{\mc W}$, which depends on $P$ and $\mc W$ only, may result the hardest to get lower bounds on in typical applications where $P$ is sparse and $\mc W$ remains small as the state space grows large. 
While Kac's formula (\cite[Lemma 21.13]{LevinPeresWilmer})
\beq\label{Kac}\sum\limits_{w\in\mc W}\sum_{v\in\mc V} \pi_wP_{wv}(\tau^v_{\mc W}+1)=1\eeq
can often be used to get upper bounds on $\tau^*_{\mc W}$ in terms of $\pi(\mc W):=\sum_{w\in\mc W}\pi_w$, lower bounds on $\tau^*_{\mc W}$ are typically harder to derive and involve finer details of $P$. 
In the last section of this paper, we will propose an analysis of $\tau^*_{\mc W}$ for networks with high local connectivity, which finds natural application when the graph associated to $P$ is a $d$-dimensional grid, and $\mc W$ is localized and its size remains bounded (or grows very slowly) as the network size grows large. Results for more general graphs, in particular, for random, locally tree-like networks will be the object of a forthcoming work.

The rest of this paper is organized as follows. 
Section \ref{sect:examples} introduces three motivating examples formalizing some of the applications mentioned at the beginning of this Introduction. 
In Section \ref{sect:mainresult}, we present our main result which is stated as Theorem \ref{theorem:main}. 
Section \ref{pagerank} discusses in detail the application of our result to the PageRank manipulation problem. Section \ref{sect:localcoonectivity} focuses on stochastic matrices whose support graph has high local connectivity and discusses lower bounds of the entrance  time $\tau^*_{\mc W}$. This allows for efficient application of Theorem \ref{theorem:main} to networks with a finite dimensional structure. Explicit examples on toroidal grid graphs are presented.


Before proceeding, let us collect here some notational conventions to be used throughout the paper. 
When referring to a graph $\mc G=(\mc V,\mc E)$, we will always use the convention that $\mc E\subseteq\mc V\times\mc V$, i.e., that its links are directed. Then, $\mc G$ undirected means that if $(u,v)\in\mc E$ then $(v,u)\in\mc E$ as well. Given $u\in\mc V$, put $\mc E_u:=\{v:(u,v)\in\mc E\}$ and let $d_u:=\!|\mc E_u|$ be the (out-) degree of node $u$. Vectors and matrices  will be considered with entries from a set $\mc V$ of finite cardinality $n:=|\mc V|$. 
The all-one column vector will be denoted by $\1$. For a matrix $A$, $A'$ will stand for its transpose and $\supp(A):=\{v:\,A_{v,\cdot}\ne0\}$ for the set of its nonzero rows. 
We refer to a probability vector as a nonnegative row vector $\mu$ such that $\mu\1=1$ and to a stochastic matrix $P$ as a nonnegative square matrix $P$ such that $P\1=\1$. 
A probability vector is said invariant for a stochastic matrix $P$  if $\mu P=\mu$. 
A stochastic matrix $P$ is said irreducible if the associated support graph $\mc G_P=(\mc V,\mc E_P)$, where $(u,v)\in\mc E_P$ if and only if $P_{uv}>0$, is strongly connected. It is a standard result that every irreducible stochastic matrix $P$ admits a unique invariant probability vector $\pi=\pi P$. 
The total variation distance between two probability vectors $\mu$ and $\pi$ is denoted by
$$||\mu-\pi||:=\frac12\sum_{v\in\mc V}|\mu_v-\pi_v|\,.$$ 
Given a stochastic matrix $P$, it is natural to consider discrete-time Markov chains $V(t)$, $t=0,1,\ldots$, with state space $\mc V$ and transition probability matrix $P$. I.e., for all $u,v\in\mc V$ and $t\ge0$,  $\P(V(t+1)=v|V(t)=u)=P_{uv}$. For $u\in\mc V$, $\P_u$ and $\E_u$ will stand for the probability and expectation conditioned on $V(0)=u$. We will also use the notation $\P_{\mu}:=\sum_{v\in\mc V}\mu_v\P_v$ for a probability vector $\mu$. We will denote the hitting time on a subset $\mc W\subseteq\mc V$ by $T_{\mc W}:=\inf\{t\ge0:\,V(t)\in\mc W\}$. It is a consequence of the Markov property that the expected hitting times $\E_u[T_{\mc W}]$ coincide with the solution $\tau^u_{\mc W}$ of equation \eqref{hittingdef}. 

\section{Three motivating applications}\label{sect:examples}
In this section we present three motivating examples formalizing some of the application problems discussed in the Introduction. Througout, $n:=|\mc V|$ will stand for the network size. 

\subsection{PageRank manipulation}\label{ex:PageRank} 



Let $\mc G=(\mc V,\mc E)$ be the directed graph describing the WWW, whose nodes $v\in\mc V$ correspond to webpages and where there is a directed link $(u,v)\in\mc E$ whenever page $u$ has a hyperlink redirecting to page $v$. 
Define a stochastic matrix $Q$ by putting $Q_{uv}=1/n$ for all $v$ if $d_u=0$, and, if $d_u\ge1$, letting $Q_{uv}=0$ if $(u,v)\notin\mc E$ and $Q_{uv}=1/d_u$ if $(u,v)\in\mc E$. Given an arbitrary probability vector $\mu$ and a parameter $\beta$ in the interval $(0,1)$, consider the equation
\be\label{pieq1}\pi=(1-\beta)\pi Q+\beta\mu\,.\ee
Since the matrix $W:=(I-(1-\beta)Q)$ is strictly diagonally dominant, hence nonsingular, equation \eqref{pieq1} admits exactly one solution $$\pi=\beta\mu W^{-1}=\beta\sum_{k\ge0}(1-\beta)^k\mu Q^k\,.$$  Observe that such vector $\pi$ turns out to be a probability vector, since each term $(1-\beta)^k\mu Q^k$ is nonnegative so that $\pi$ is as well, and $\mu Q^k\1=\mu\1=1$ for every $k\ge0$, so that $\pi\1=\beta\sum_{k\ge0}(1-\beta)^k=1$. The vector $\pi$ is known as the PageRank vector and was first introduced by Brin and Page \cite{BrinPage:98} to measure the relative importance of webpages. In the original PageRank version, $\mu=n^{-1}\1$ is chosen as the uniform distribution over the set of webpages, while typical values of $\beta$ used in practice are about $0.15$. More general choices of the probability vector $\mu$ lead to the definition of the personalized PageRank \cite{Haveliwala:2003}, which  is used in context-sensitive searches. 

Consider now a (relatively small) set of webpages $\mc W\subseteq\mc V$, and assume that the set  $\bigcup_{w\in\mc W}\mc E_w$ of hyperlinks originated from these webpages can be modified arbitrarily in order to change $\pi$. Let $\tilde{\mc G}=(\mc V,\tilde{\mc  E})$ be the modified WWW graph, $\tilde Q$ the corresponding stochastic matrix, and $\tilde\pi$ the corresponding modified PageRank vector solving the equation
\be\label{pieq3}\tilde\pi=(1-\beta)\tilde\pi\tilde Q+\beta\mu\,.\ee

We now give a different characterization of the PageRank vector and reformulate the perturbation problem.
First, we introduce the stochastic matrix
$$P:=(1-\beta)Q+\beta\1\mu\,.$$ 
We claim that $P$ has a unique invariant probability vector and that it coincides with the PageRank vector $\pi$.
To see this equivalence,  notice that, 
if $\nu$ is any row vector such that $\nu\1=1$, we have that 
$$\nu P=(1-\beta)\nu Q+\beta\nu\1\mu=(1-\beta)\nu Q+\beta\mu\,,$$
so that $\nu=\nu P$  if and only if $\nu$ coincides with the solution $\pi$ of \eqref{pieq1}. An analogous argument shows that the modified PageRank vector $\tilde\pi$ coincides with the unique invariant probability vector of the stochastic matrix 
$$\tilde P:=(1-\beta)\tilde Q+\beta\1\mu\,.$$

A standard result \cite[Proposition 4.2]{LevinPeresWilmer} allows one to write the total variation distance between $\pi$ and $\tilde\pi$ as  
\be\label{TVidentity}||\tilde\pi-\pi||=\max_{\mc U\subseteq\mc V}\l\{\tilde\pi(\mc U)-\pi(\mc U)\r\}\,.\ee
Hence, estimating the impact that an arbitrary change of the hyperlinks from a subset $\mc W$ of webpages has on the aggregate PageRank of an arbitrary subset $\mc U$ of webpages boils down to bounding the total variation distance between the invariant probability vectors $\pi$ and $\tilde\pi$ of the stochastic matrices $P$ and $\tilde P$, respectively. Observe that, since the matrices $Q$ and $\tilde Q$ differ only on the rows indexed by elements of $\mc W$, so do $P$ and $\tilde P$. 
In Example \ref{ex:PageRank} of Section~\ref{sect:mainresult} we will prove an upper bound on $||\tilde\pi-\pi||$ depending only on the size of $\mc W$ (as measured by $\pi$ and $\mu$), and on the value of the parameter $\beta\in(0,1)$.

\subsection{Faulty communication links in distributed averaging algorithms} Consider a sensor network described as a connected undirected  graph $\mc G=(\mc V, \mc E)$, whose nodes and links represent sensors and two-way communication links, respectively. 
Assume that each sensor $v$ initially measures a scalar value $y_v$ and the goal is to design a distributed algorithm for the computation of the arithmetic average 
$$\ov y:=\frac1n\sum_{v\in\mc V}y_v\,.$$ 

A possible solution \cite{OlshevskyTsitsiklis:2011} is as follows. Let $d\in\R^{\mc V}$ be the degree vector in $\mc G$. Initialize the state of every sensor $v\in\mc V$ as
\be\label{averaging0}x_v(0)=\frac{y_v}{d_v}\,,\qquad z_v(0)=\frac1{d_v}\,.\ee
Then, at every time instant $t=0,1,\ldots$, let every sensor $v\in\mc V$ update its state according to the recursion 
\be\label{averaging1} [x_v(t+1),z_v(t+1)]\!=\!\frac12[x_v(t),z_v(t)]+\frac1{2d_v}\!\!\sum_{\substack{u\in\mc V:\,\\(v,u)\in\mc E}}\!\![x_u(t),z_u(t)]\,.\ee
What makes the above iteration particularly appealing in large-scale network applications is the fact that it requires sensors to exchange information with their neighbors in $\mc G$ only, and that each sensor $v$ needs to know only its degree $d_v$ and initial measurement $y_v$, with no need for global knowledge about  the network structure or size. 

In order to analyze the algorithm let us rewrite (\ref{averaging0}) and (\ref{averaging1}) in matrix notation. 
Let $P$ be the stochastic matrix associated to the lazy random walk on $\mc G$, i.e., $P=(I+Q)/2$, where $I$ denotes the identity matrix and $Q_{uv}=1/d_u$ if $(u,v)\in\mc E$ and $Q_{uv}=0$ if $(u,v)\notin\mc E$.  
Let 
\be\label{averaging0.1}x(0)=\frac yd\,,\qquad z(0)=\frac\1d\ee (where division between two vectors is meant componentwise) and consider the iteration 
\be\label{averaging1.1}x(t+1)=Px(t)\,,\qquad z(t+1)=Pz(t)\,.\ee
Observe that the unique invariant probability vector $\pi$ of the matrix $P$ is given by \be\label{explicitpi}\pi_v=\frac{d_v}{n\ov d}\,,\qquad v\in\mc V\,,\ee where  
$$ \ov d:=\frac1n\sum_{v\in\mc V} d_v$$ is the average degree. 
Moreover, irreducibility and acyclicity of $P$ (implied by $P_{uu}>0$ for all $u$) imply that 
$$x(t)=P^t\frac yd\stackrel{t\to\infty}{\longrightarrow}\1\pi\frac yd=\1\frac{\ov y}{\ov d}\,,$$
$$ z(t)=P^t\frac{\1}{d}\stackrel{t\to\infty}{\longrightarrow}\1\pi\frac{\1}d=\1\frac1{\ov d}\,,$$ 
so that 
$$\frac{x_v(t)}{z_v(t)}\stackrel{t\to\infty}{\longrightarrow}\ov y\,,\qquad \forall v\in\mc V\,.$$
Therefore, the iterative distributed algorithm defined by (\ref{averaging0.1})-(\ref{averaging1.1}) effectively computes the average $\ov y$ of the vector $y$. The example can be generalized to those weighted graphs whose nodes all have in-degree equal to the out-degree (hence, in particular, undirected weighted graphs). Indeed, for these graphs, the invariant probability vector $\pi$ of the associated stochastic matrix $P$ admits the explicit form \eqref{explicitpi}.

Now, let $\mc F\subseteq\mc E$ be a subset of directed communication  links which stop working and $\tilde{\mc G}:=(\mc V,\tilde{\mc E})$, where $\tilde{\mc E}:=\mc E\setminus\mc F$, be the directed graph obtained from $\mc G$ by removing such links. Let $\tilde d$ be the vector of in-degrees in $\tilde{\mc G}$ and define $\tilde P=(I+\tilde Q)/2$, where $\tilde Q$ is a stochastic matrix with $\tilde Q_{uv}=1/\tilde d_u$ if $(v,u)\in\tilde{\mc E}$ and $\tilde Q_{uv}=0$ otherwise. 
Consider the following recursion, analogous to (\ref{averaging0.1}) and (\ref{averaging1.1}), with $d$ and $\mc E$ replaced by $\tilde d$ and $\tilde{\mc E}$, respectively:  
\be\label{averaging2}\tilde x(0)= \frac y{\tilde d}\,,\qquad\tilde z(0)=\frac{\1}{\tilde d}\,,\ee  
\be\label{averaging3}\tilde x(t+1)=\tilde Px(t)\, \qquad\tilde z(t+1)=\tilde P\tilde z(t)\,.\ee Then, provided that $\tilde{\mc G}$ remains strongly connected, an argument as the one before shows that 
$$\frac{\tilde x_v(t)}{\tilde z_v(t)}\stackrel{t\to\infty}{\longrightarrow}\tilde y\,,\qquad \forall v\in\mc V\,,$$
where 
$$\tilde y=\frac{\tilde\pi(y/\tilde d)}{\tilde\pi(\1/\tilde d)}$$ 
and
$\tilde\pi$ is the unique invariant probability vector of $\tilde P$.
In other words, the perturbed dynamics \eqref{averaging2}--\eqref{averaging3} achieve consensus on a perturbed value $\tilde y$.

We are now going to show that the absolute error $|\tilde y-y|$ can be upper bounded in terms of the total variation $||\tilde\pi-\pi||$ and the fraction $|\mc F|/|\mc E|$ of failed communication links. To see this, first we express the perturbed consensus value as
$$\tilde y=\frac{\tilde\pi y/\tilde d}{\tilde\pi\1/\tilde d}=\frac{\ov y+\eps_1+\eps_2}{1+\eps_3+\eps_4}\,,$$
where 
$$\eps_1:=\frac1n\sum_{v\in\mc V}\l(\frac{d_v}{\tilde d_v}-1\r)y_v\,,\qquad\eps_2:=\ov d\sum_{v\in\mc V}(\tilde\pi_v-\pi_v)\frac{y_v}{\tilde d_v}\,,$$
$$\eps_3:=\frac1n\sum_{v\in\mc V}\l(\frac{d_v}{\tilde d_v}-1\r)\,,\qquad\eps_4:=\ov d\sum_{v\in\mc V}(\tilde\pi_v-\pi_v)\frac{1}{\tilde d_v}\,.$$
Now, using the facts that $\tilde d_v\ge1$ for all $v$ (since $\tilde{\mc G}$ is connected) and that $|\mc E|=\sum_{v\in\mc V}d_v=n\ov d$, one gets that 
$$|\eps_1|\le \ov d\frac{|\mc F|}{|\mc E|}||y||_{\infty}\,,\quad |\eps_2|\le \ov d||y||_{\infty}||\tilde\pi-\pi||_{}\,,$$
$$|\eps_3|\le \ov d\frac{|\mc F|}{|\mc E|}\,,\quad |\eps_4|\le \ov d||\tilde\pi-\pi||_{}\,,\quad |\tilde y|\le||y||_{\infty}\,.$$
It follows that 
$$\ba{rcl}|\tilde y-\ov y|&=&\l|\eps_3\tilde y+\eps_4\tilde y-\eps_1-\eps_2\r|\\[5pt]
&\le&
|\tilde y|(|\eps_3|+|\eps_4|)+|\eps_1|+|\eps_2|
\\[5pt]
&\le&
2\ov d||y||_{\infty}\l(|\mc F|/|\mc E|+||\tilde\pi-\pi||\r)
\,,\ea$$
so that 
\be\label{boundy-y}
\frac{|\tilde y-\ov y|}{||y||_{\infty}}\le2\ov d\l(\frac{|\mc F|}{|\mc E|}+||\tilde\pi-\pi||\r)
\,.\ee
Formula \eqref{boundy-y} shows that, provided an upper bound on the average degree $\ov d$, in order to guarantee that the value $\tilde y$ computed by the distributed averaging algorithm on the perturbed graph $\tilde{\mc G}$ is close to the average $\ov y$ of the sensors' measurements, it is sufficient that both the fraction $|\mc F|/|\mc E|$ of failed communication links and the total variation distance $||\tilde\pi-\pi||$ are small. 

\subsection{Voter model with influential agents} \label{ex:minoritiesnew}
Let $\mc G=(\mc V,\mc E)$ be a connected undirected graph (with no self-loops). Nodes are to be interpreted as agents possessing a binary opinion. Opinions change in time as a consequence of pairwise interactions in the network. Precisely, for $u\in\mc V$ and $t=0,1,\ldots$, let $X_u(t)\in\{0,1\}$ be the opinion of agent $u$ at time $t$. At every time $t=0,1,\ldots$, a single directed link $(u,v)$ is activated, chosen uniformly at random from $\mc E$, and its tail node $u$ updates its state $X_u(t)$ by copying the head node $v$'s current state $X_v(t)$. 
By assembling all the agents' opinions in a vector $X(t)\in\{0,1\}^{\mc V}$ we obtain that $X(t)$ is a Markov chain whose transitions can be compactly described as follows. For $u\ne v\in\mc V$, let $E^{(u,v)}\in\R^{\mc V\times\mc V}$ have all entries equal to zero but for $E^{(u,v)}_{u,v}=-E^{(u,v)}_{u,u}=1$. 
Then, given $X(t)$, we have that $$X(t+1)=(I+E^{(u,v)})X(t)$$ with probability $1/|\mc E|$, for all $(u,v)\in\mc E$. 
This is an instance of the voter model \cite{Liggettbook,Liggettbook2,ConnelyWelsh82,Cox89}. In a social network interpretation, this may be thought of modeling a society where every pair of individuals whose corresponding nodes are neighbors in $\mc G$ have the same chance to influence each other. It is a standard result that, with probability one, this dynamics achieves consensus in some finite time. More precisely, there exists some random consensus time $T$, which is finite with probability one, and a random consensus value $Y\in\{0,1\}$, such that 
\be\label{consensusfinitetime}X_v(t)=Y\,,\qquad v\in\mc V\,,\  t\ge T\,.\ee
The main asymptotic quantity of interest is the probability distribution of the consensus value $Y$ conditioned to the initial condition $X(0)$. Specifically, we define
$$y:=\P(Y=1|X(0))\,.$$


Now, let us consider the following variant to the model. Consider a directed subgraph $\tilde{\mc G}=(\mc V,\tilde{\mc E})$, where $\tilde{\mc E}=\mc E\setminus\mc F$ is  obtained from $\mc E$ by removing a subset  $\mc F\subseteq\mc E$ of directed links. 
We assume that $\tilde{\mc G}$ remains strongly connected. 
Consider the Markov chain $\tilde X(t)$ over $\{0,1\}^{\mc V}$ such that, given $\tilde X(t)$, 
$$\tilde X(t+1)=(I+E^{(u,v)})\tilde X(t)$$ with probability $|\mc E|^{-1}$, for all $(u,v)\in\tilde{\mc E}$, and $\tilde X(t+1)=\tilde X(t)$ with probability $|\mc F|/|\mc E|$. The social network interpretation is that $$\mc W:=\{u:\,(u,v)\in\mc F\text{ for some }v\}$$ is a set of influential individuals, whose interactions with some of their neighbors in $\mc G$ are asymmetric, as they influence such neighbors without being influenced in turn from them.  A similar model is discussed in \cite{influence} in the framework of opinion dynamics over continuous space. Observe that, analogously to the voter model on $\mc G$, strong connectivity of the graph $\tilde{\mc G}$ implies that, with probability one, the process $\tilde X(t)$ achieves a consensus in finite time on a binary random variable $\tilde Y$. We can similarly define the conditional probability
$$\tilde y:=\P(\tilde Y=1|\tilde X(0))\,.$$

The absolute difference $|\tilde y-y|$ measures the effect of the influential individuals in the final consensus value.
We now give a different characterization for $y$ and $\tilde y$ in terms of invariant probability vectors of suitably defined stochastic matrices and propose a characterization of $|\tilde y-y|$ in terms of their total variation difference.

Let us define the stochastic matrix
$$P:=I+\frac{1}{|\mc E|}\sum_{(u,v)\in\mc E}E^{(u,v)}\,.$$
Then, $$\E[X(t+1)|X(t)]=PX(t)\,,\qquad t=0,1,\ldots\,,$$ so that an inductive argument proves that 
\be\label{Xtplus}\E[X(t)|X(0)]=P^tX(0)\,\qquad t\ge0\,.\ee
Since $\mc G$ is connected and undirected, $P$ is irreducible and symmetric, so that its unique invariant probability vector is the uniform one 
$$\pi=\frac1n\1'\,.$$
It then follows from \eqref{Xtplus} that, for all $t\ge0$,  
\be\label{Xtplus1}
\ba{rcl}
\ds\E\l[\frac1n\sum_{v\in\mc V} X_v(t)|X(0)\r]
&=&\ds\pi\E[X(t)|X(0)]
\\&=&\ds\pi P^tX(0)\\[5pt]
&=&\ds\pi X(0)\\[5pt]
&=&\ds\frac1n\sum_{v\in\mc V}X_v(0)\,,\ea\ee
a property that is sometimes referred to as 
as conservation of the average magnetization \cite{Suchecki.etal:2005} in the statistical physics jargon. 
Finally, it follows from \eqref{consensusfinitetime} and \eqref{Xtplus1} that 
$$
y=\E[Y|X(0)]=\lim_{t\to+\infty}\frac1n\sum_{v\in\mc V}\E[X_v(t)|X(0)]=\frac1n\sum_{v\in\mc V}X_v(0)\,.
$$
Similarly,
$$\tilde y=\tilde\pi\tilde X(0)\,,$$
where $\tilde\pi=\tilde\pi \tilde P$ is the unique invariant probability vector of the stochastic matrix
$$\tilde P:=I+\frac1{|\mc E|}\sum_{(u,v)\in\tilde{\mc E}}E^{(u,v)}\,.$$ 
Clearly, if the initial conditions of the two processes coincide, i.e., if $\tilde X(0)=X(0)$, then 
$$|\tilde y-y|\le||\tilde\pi-\pi||\,.$$
In fact, while the inequality above is valid for every initial state value $\tilde X(0)=X(0)\in\{0,1\}^{\mc V}$, the identity \eqref{TVidentity} implies that such inequality is tight in the sense that there exists one value $x\in\{0,1\}^{\mc V}$ (the one with $x_u=1$ for $u\in\mc U$ and $x_v=0$ for $v\in\mc V\setminus\mc U$, where $\mc U$ is such that $||\tilde\pi-\pi||=\tilde\pi(\mc U)-\pi(\mc U)$) such that, if $\tilde X(0)=X(0)=x$, then $|\tilde y-y|=||\tilde\pi-\pi||$.

It follows that the problem of estimating the difference between the probability vector of the eventual consensus value for the voter model on $\mc G$ and $\tilde{\mc G}$ is equivalent to the one of estimating the total variation distance between the invariant probability vectors of the stochastic matrices $P$ and $\tilde P$, respectively.

\section{Perturbation results}\label{sect:mainresult}
Let $P$ be an irreducible stochastic matrix on the finite state space $\mc V$ 
and let $\pi=\pi P$ be its unique invariant probability vector. Let $\tilde P$ be another stochastic matrix (not necessarily irreducible) on the same state space $\mc V$, to be interpreted as a perturbation of $P$, and let $\tilde\pi=\tilde\pi\tilde P$ be an invariant probability vector of $\tilde P$ (not necessarily the unique one). 



The following result provides an upper bound on the total variation distance between $\pi$ and $\tilde \pi$.
\begin{figure}\begin{center}
\includegraphics[height=4.5cm,width=6.0cm]{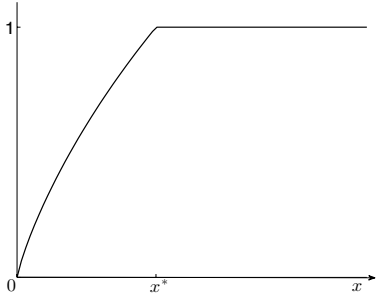}
\caption{\label{fig:thetaplot}Graph of the function $\Psi(x)$ defined in (\ref{eq:thetadef}).}
\end{center}
\end{figure}
It is stated in terms of the function $\Psi:[0,+\infty)\to[0,1]$
\be\label{eq:thetadef}
\Psi(x):=\l\{\ba{lcl}x\ln\l(e^2/x\r)&\ & x\le x^*\\ 1&\ &x> x^*\,,\ea\r.\ee
where $x^*=0.31784\ldots$ is the smallest positive solution of $x\ln (e^2/x)=1$. (The graph of $\Psi(\,\cdot\,)$ is plotted in Figure \ref{fig:thetaplot}.)
\begin{lemma}\label{lemma:TVestimate} 
Let $P$ and $\tilde P$ be stochastic matrices on a finite set $\mc V$. 
Let $P$ be irreducible with invariant probability vector $\pi$ and mixing time $\tmix$ as defined in \eqref{eq:mixingdef}, and $\tilde\pi$ be an invariant probability vector for $\tilde P$. Then, 
$$||\tilde\pi-\pi||_{}\le\Psi(\tmix\cdot\,\tilde\pi(\mc W))\,,$$
for all $\mc W\subseteq\mc V$ such that $\mc W\supseteq\supp(P-\tilde P)$.
\end{lemma}
\proof 
Let $V(t)$ and  $\tilde V(t)$ be two Markov chains on $\mc V$ which start and move together with transition probabilities $P_{uv}$ until the first time $T_{\mc W}=\tilde T_{\mc W}$ they hit $\mc W$, and move independently with transition probabilities $P_{uv}$ and $\tilde P_{uv}$, respectively, ever after. Since $P$ and $\tilde P$ coincide on $\mc V\setminus\mc W$, one has that the marginal transition probability matrices of $V(t)$ and $\tilde V(t)$ coincide with $P$ and $\tilde P$, respectively. 
Then, for all $\mc A\subseteq\mc V$, and $t\ge0$, one has that
$$\ba{rcl}\!\!\tilde\pi(\mc A)\!\!&=&\P_{\tilde\pi}(\tilde V(t)\in\mc A)\\[10pt]
&=&\P_{\tilde\pi}(\tilde V(t)\in\mc A,\tilde T_{\mc W}\ge t)+\P_{\tilde\pi}(\tilde V(t)\in\mc A,\tilde T_{\mc W}<t)\\[10pt]
&=&\P_{\tilde\pi}(V(t)\in\mc A,\tilde T_{\mc W}\ge t)+\P_{\tilde\pi}(\tilde V(t)\in\mc A,\tilde T_{\mc W}<t)\\[10pt]
&\le&\P_{\tilde\pi}(V(t)\in\mc A)+\P_{\tilde\pi}(\tilde T_{\mc W}< t)\\[10pt]
&\le&\pi(\mc A)+\exp(-\lfloor t/\tmix\rfloor)+ t\tilde\pi(\mc W)\,,
\ea$$
where the first identity uses the invariance of $\tilde\pi$, and the last inequality follows from 
$$
\ba{rcl}\P_{\tilde\pi}(V(t)\in\mc A)-\pi(\mc A)&=&
[\tilde\pi P^t](\mc A)-\pi(\mc A)\\[5pt]&\le&
||\tilde\pi P^t-\pi||_{}\\[5pt]&\le&
\exp(-\lfloor t/\tmix\rfloor)\,,\ea$$ which is a consequence of the representation \eqref{TVidentity} of the total variation distance and of the submultiplicativity  property of the maximal total variation distance, see, e.g., they discussion following formula (4.31) in \cite{LevinPeresWilmer}), and  the bound  $$\P_{\tilde\pi}(\tilde T_{\mc W}< t)\le\sum_{i=0}^{t-1}\P_{\tilde\pi}(\tilde V(i)\in\mc W)=t\tilde\pi(\mc W)\,,$$ 
which is implied by the union bound and, again, invariance of $\tilde\pi$ for $\tilde P$.
Therefore, using the characterization \eqref{TVidentity} of the total variation distance, one gets that
$$||\tilde\pi-\pi||_{}=\max_{\mc A\subseteq\mc V}\l\{\tilde\pi(\mc A)-\pi(\mc A)\r\}\le\exp(-\lfloor t/\tmix\rfloor)+ t\tilde\pi(\mc W)\,, $$
for all $t\ge0$. 
The claim now follows by choosing $$t=\max\l\{\l\lfloor\tmix\log\frac{e}{\tmix\cdot\tilde\pi(\mc W)}\r\rfloor,0\r\}\,,$$ such a choice being suggested by the minimization of the function $$x\mapsto\exp(-x/\tmix-1)+ x\tilde\pi(\mc W)$$ over continuous nonnegative values of $x$. 
\qed

Lemma \ref{lemma:TVestimate} shows that it is sufficient to have an upper bound on the product $\tmix\cdot\tilde\pi(\mc W)$ in order to obtain an upper bound on the total variation distance $||\tilde\pi-\pi||_{}$. 
In particular, assuming that an upper bound on the mixing time $\tmix$ is available, e.g., from an estimate of the conductance of $P$, one is left with estimating $\tilde\pi(\mc W)$. Observe that $\tilde\pi(\mc W)$ is typically unknown in the applications. 
Below, we derive an upper bound on $\tilde\pi(\mc W)$ in terms of the entrance time $\tau^*_{\mc W}$ and of the exit probability $\tilde\gamma_{\mc W}$, defined in \eqref{eq:tauWdef} and \eqref{def:exitprobability}, respectively. These two quantities can be given the following probabilistic interpretation. Consider a Markov chain $\tilde V(t)$ on $\mc V$ with transition probability matrix $\tilde P$, and let $$\tilde T_{\mc W}:=\inf\{t\ge0:\,\tilde V(t)\in\mc W\}$$ and $$\tilde T_{\mc V\setminus\mc W}:=\inf\{t\ge0:\,\tilde V(t)\in\mc V\setminus\mc W\}$$ be, respectively, the hitting time on, and the exit time from, the set $\mc W$. 
Then, since $P$ and $\tilde P$ coincide outside $\mc W$, one has that the expected hitting times 
satisfy 
\be\label{identityhittingtimes}\E_v[\tilde T_{\mc W}]=\tau^v_{\mc W}=\E_v[T_{\mc W}]\,,\qquad v\in\mc V\,.\ee 
In fact, 
the entrance time $\tau^*_{\mc W}=\min\{\tau^u_{\mc W}:\,v\in\mc V\setminus\mc W\}$ only depends on the choice of the subset $\mc W\supseteq\supp(\tilde P-P)$ and on the original matrix $P$ (in particular, on the rows of $P$ indexed by $v\notin\mc W$), but not on finer details of the perturbation $\tilde P-P$. On the other hand, for every $w\in\mc W$ and $k\ge1$ one has that 
$$\phi_w(k):=\P_w(\tilde T_{\mc V\setminus\mc W}=k)=\sum_{\substack{\xi_0=w,\xi_k\in\mc V\setminus\mc W\\\xi_1,\ldots, \xi_{k-1}\in\mc W}}\ \prod_{1\le l\le k}\tilde P_{\xi_{l-1}\xi_l}\,,$$
so that the exit probability defined in \eqref{def:exitprobability} satisfies 
\be\label{eq:kappaWdef}\tilde\gamma_{\mc W}=\sup_{t\ge1}\,\min_{\substack{w\in\mc W:\\[1pt]\tilde\pi_w>0}}\,\frac1t\sum_{1\le k\le t}\phi_w(k)
=\sup_{t\ge1}\ \min_{\substack{w\in\mc W:\\[1pt]\tilde\pi_w>0}}\!\frac{\P_w(\tilde T_{\mc V\setminus\mc W}\le t\big)}t
.\ee
Notice that the exit probability $\tilde\gamma_{\mc W}$ depends only on those rows of the perturbed matrix $\tilde P$ whose indices lie in $\mc W$ (because so does the distribution  of $\tilde T_{\mc V\setminus\mc W}$) and, when $\tilde P$ is not irreducible, on the choice of the invariant measure $\tilde\pi$. In particular, one has that $\tilde\gamma_{\mc W}=0$ if and only if $\mc V\setminus\mc W$ is not accessible under $\tilde P$ from some state $w\in\mc W$ such that $\tilde\pi_w>0$. 
We are now in a position to prove the following result.
\begin{lemma}\label{lemma:W}
Let $\tilde P$ be a stochastic matrix on a finite set $\mc V$, and $\tilde\pi=\tilde\pi\tilde P$ an invariant probability measure. 
Then,
\begin{equation}\label{estimbott}\tilde\pi(\mc W)\le\frac{1}{\tilde\gamma_{\mc W}\cdot\tau^*_{\mc W}}\,,\end{equation}
for all $\mc W\subseteq\mc V\,.$
\end{lemma}
\proof  
Observe that, for $k\ge1$ and $w\in\mc W$,
\be\label{usefulLemma2}\sum_{v\in\mc V}\tilde P_{wv}\tau^v_{\mc W}=
\sum_{u\in\mc V\setminus\mc W}\tilde P_{wu}\tau^u_{\mc W}
\ge\tau^*_{\mc W}\phi_w(1)\,.\ee
Then, it follows from Kac's formula \eqref{Kac} applied to $\tilde P$ and $\tilde\pi$, the identity \eqref{identityhittingtimes}, and the inequality \eqref{usefulLemma2}, that 
\be\label{estimate1}\ba{rcl}\ds\frac1{\tilde\pi(\mc W)}-1
&=&\ds\frac1{\tilde\pi(\mc W)}\sum_{w\in\mc W}\sum_{v\in\mc V}\tilde\pi_w\tilde P_{wv}\tau^v_{\mc W}\\[15pt]
&\ge&\ds\frac{\tau^*_{\mc W}}{\tilde\pi(\mc W)}\sum_{w\in\mc W}\tilde\pi_w\phi_w(1)\,.\ea\ee
%
%
%
Now, observe that 
$$\sum_{w'\in\mc W}\tilde\pi_{w'}\tilde P_{w'w}\leq\sum_{v\in\mc V}\tilde\pi_{v}\tilde P_{vw}= \tilde\pi_{w}\,, \qquad w\in\mc W\,.$$
Then, for all $k\ge1$, one gets that
$$\ba{rcl}\ds\sum_{w'\in\mc W}\tilde\pi_{w'}\phi_{w'}(k+1)&=&
\ds\sum_{w'\in\mc W}\sum_{w\in\mc W}\tilde\pi_{w'}\tilde P_{w'w}\phi_{w}(k)\\[15pt]
&\le&\ds\sum_{w\in\mc W}\tilde\pi_{w}\phi_{w}(k)\,.\ea$$
It follows that, for all $t>0$, 
\be\label{estimate2}\ds\sum_{w\in\mc W}\tilde\pi_w\phi_w(1)\ge
\sum_{w\in\mc W}\tilde\pi_w\,\cdot\,\ds\frac1t\sum_{1\le k\le t}\phi_w(k)\,.\ee
Then,  (\ref{estimate1}) and (\ref{estimate2}) imply that 
$$\ba{rcl}\ds\frac1{\tilde\pi(\mc W)}
&\ge& \ds\tau^*_{\mc W}\sum_{w\in\mc W}\frac{\tilde\pi_w}{\tilde\pi(\mc W)}\,\cdot\, \frac{1}{t}\sum\limits_{k=1}^t \phi_w(k)\\[5pt]
&\geq&\ds\tau^*_{\mc W}\min_{\substack{w\in\mc W:\\[1pt]\tilde\pi_w>0}}
\frac{1}{t}\sum\limits_{k=1}^t \phi_w(k)\,.\ea
$$
Since $t\ge1$ is arbitrary, the inequality above implies that 
$$\frac{1}{\tilde\pi(\mc W)}\ge\tau^*_{\mc W}\cdot\tilde\gamma_{\mc W}\,,$$
thus proving the claim.
\qed

Lemmas \ref{lemma:TVestimate} and \ref{lemma:W} immediately imply the following result:
\begin{theorem} \label{theorem:main}
Let $P$ and $\tilde P$ be stochastic matrices on a finite set $\mc V$. 
Let $P$ be irreducible with invariant probability measure $\pi$ and mixing time $\tmix$, and $\tilde\pi$ be an invariant probability measure for $\tilde P$. Then, 
$$||\tilde\pi-\pi||_{}\le\Psi\left(\frac{\tmix}{\tilde\gamma_{\mc W}\cdot\tau_{\mc W}^*}\right)\,,$$
for all $\mc W\subseteq\mc V$ such that $\supp(\tilde P-P)\subseteq\mc W$.
\end{theorem}

Theorem \ref{theorem:main} implies that, in order for the total variation distance $||\tilde\pi-\pi||$ to be small, it is sufficient that, for some set $\mc W\supseteq\supp(\tilde P-P)$, the ratio
$$\frac{\tmix}{\tilde\gamma_{\mc W}\cdot\tau_{\mc W}^*}$$ is small.  

While the term $1/\tau_{\mc W}^*$ is a measure of `how large' the set $\mc W$ is, hence it is expected to play a central role in estimating $||\tilde\pi-\pi||$, one might wonder what the roles of the mixing time $\tmix$ and of the exit probability $\tilde\gamma_{\mc W}$ are. 
The following two simple examples show that having control of each of the terms $\tilde\gamma_{\mc W}$ and $\tmix$ is indeed necessary in order to bound the total variation distance $||\tilde\pi-\pi||$.


\begin{figure}\begin{center}
{\includegraphics[height=4.2cm,width=9.0cm]{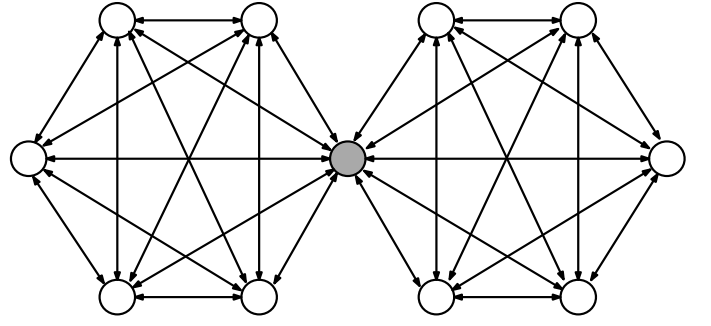}}
\caption{\label{figglued}The graph of Example \ref{example:gluedcomplete}, for $m=5$. The perturbation set $\mc W=\{0\}$ is shaded in gray.}
\end{center}
\end{figure}

\begin{example}
For an integer $n\ge 2$, consider the stochastic matrix $P$ of size $n\times n$, with all entries equal to $1/n$. Perturb it in a single row $w$ by putting, for some $\alpha\in(0,1-1/n)$, 
$$\tilde P_{ww}=1-\alpha\,,\qquad\tilde P_{wv}=\alpha/(n-1)\,,\quad v\ne w\,.$$ 
Then, $\tmix=1$, $\tau^*_{\mc W}=n$, and $\tilde\gamma_{\mc W} =\alpha$. Then, it follows from Theorem \ref{theorem:main} that $\alpha n\to\infty$
is a sufficient condition for $||\tilde\pi-\pi||\to0$ as $n$ grows large. 
On the other hand, it is easily verified that $\pi_v=1/n$ for all $v$, while $\tilde\pi_w=1/(n\alpha+1)$, and $\tilde\pi_v=n\alpha/((n-1)(n\alpha+1))$, for all $v\ne w$. Hence, 
$$||\tilde\pi-\pi||_{}=\frac{1-\alpha-1/n}{n\alpha+1}$$ 
which shows that $\alpha n\to\infty$ is indeed also a necessary condition for $||\tilde\pi-\pi||\to0$ as $n$ grows large. %
\end{example}


\begin{example}\label{example:gluedcomplete}
For a positive integer $m$, define the stochastic matrix $P$ on the set $\mc V:=\{-m,-m+1,\ldots,m-1,m\}$ by putting $P_{uv}=1/m$ if $u\ne v$ and $u\cdot v\ge0$, $P_{uv}=0$ if $u\cdot v<0$ or $u=v$, and $P_{0v}=1/(2m)$ for all $v\ne0$. 
Such $P$ can be interpreted as the transition probability matrix associated to the random walk on the graph of Figure \ref{figglued}. 
Then, one has that 
$$\pi_0=\frac1{m+1}\,,\qquad\pi_v=\frac1{2m+2}\,,\quad v\ne0\,.$$ 
Now, for some $0<\alpha<1/2$, perturb $P$ on $\mc W=\{0\}$ by putting 
$$\tilde P_{00}=0\,,\qquad\tilde P_{0v}=\frac{1/2+\alpha\mathrm{sgn}(v)}m\,,\quad v\ne0\,.$$ Straightforward computations show that 
$$\tau^*_{\mc W}=m\,,\qquad \tilde\gamma_{\mc W}=1\,.$$ 
On the other hand, the bottleneck bound \cite[Theorem 7.3]{LevinPeresWilmer} implies that 
$$\tmix\ge\frac1{4\pi_0}\ge\frac m2\,,$$ so that 
Theorem \ref{theorem:main} is useless as it only provides the trivial conclusion that $||\tilde\pi-\pi||\le1$. 
In fact, observe that $$\tilde\pi_v-\pi_v=\frac\alpha{m+1}\mathrm{sgn}(v)\,,\qquad v\in\mc V\,,$$ so that 
$$||\tilde\pi-\pi||=m\cdot\frac{\alpha}{m+1}$$ is arbitrarily close to $\alpha$ for large $m$.
Hence, $||\tilde\pi-\pi||$ does not vanish as $m$ grows large, unless $\alpha$ itself does so. 
The intuitive explanation is that, while the perturbation is concentrated on a single node, $w=0$, that is assigned a relatively small weight $\pi_0$ by the invariant probability vector $\pi$, such node lies along every path connecting two nodes $u$ and $v$ of opposite sign. This is reflected in the large mixing time $\tmix$.

\end{example}

\begin{figure}\begin{center}
\hspace{-0cm}
\includegraphics[height=4.1cm,width=4.1cm]{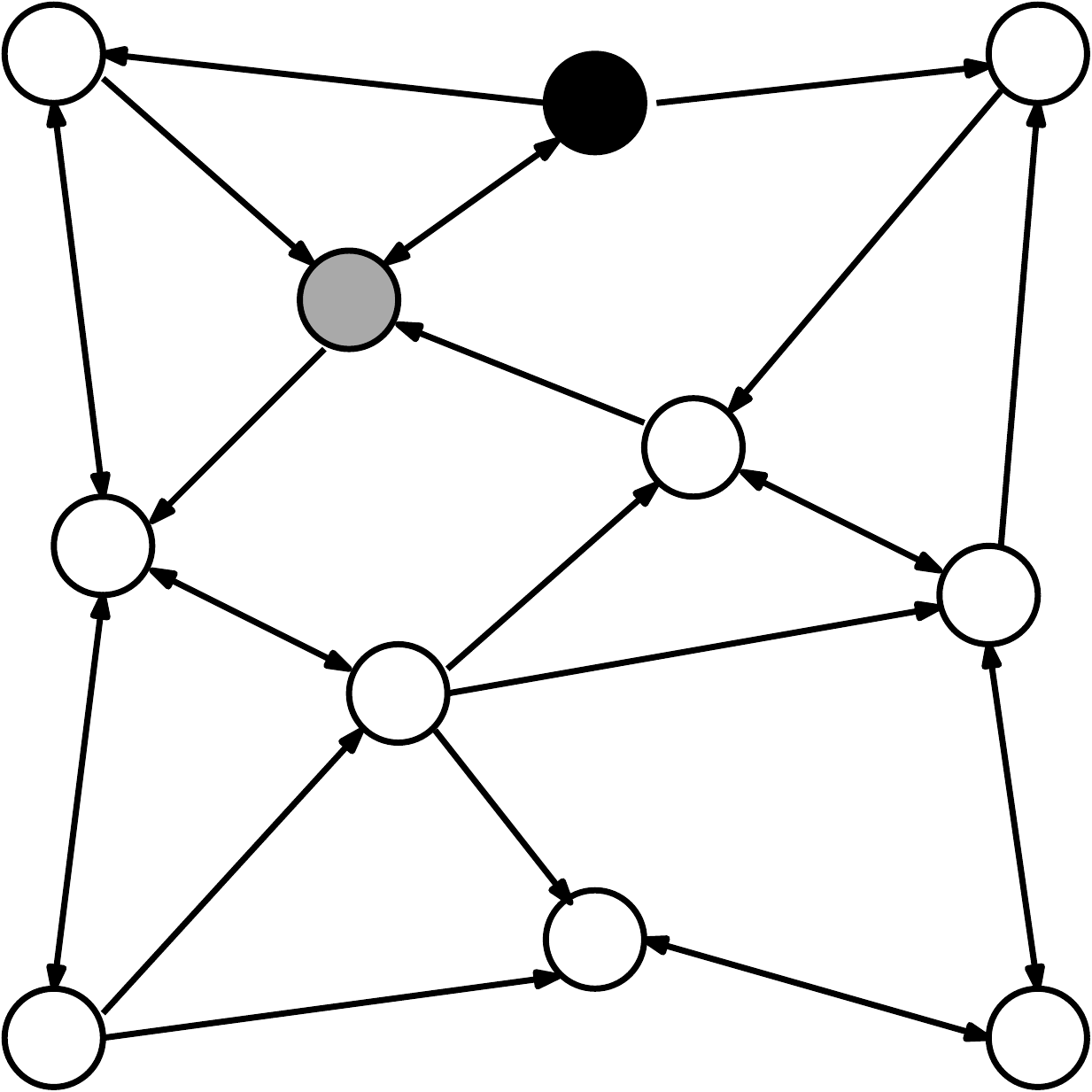}\hspace{.4cm}
\includegraphics[height=4.1cm,width=4.1cm]{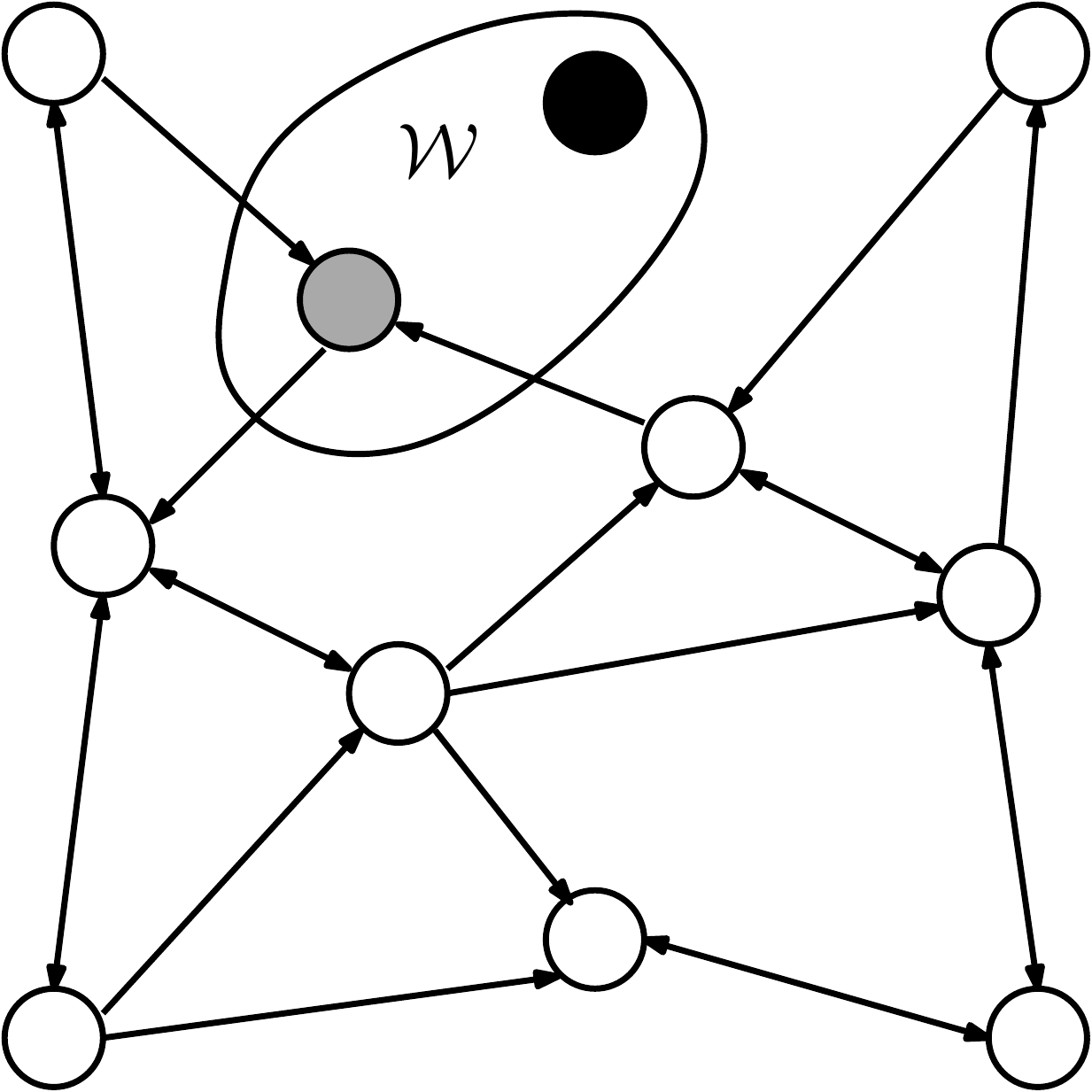}
\end{center}
\caption{\label{figdisconnected} 
On the left, the support graph $\mc G_{P}$ of an irreducible stochastic matrix $P$. On the right, the support graph $\mc G_{\tilde P}$ of a perturbed matrix $\tilde P$ obtained as in Example \ref{example:disconnected}, with node $u$ colored in black and the rest of the perturbation set $\mc W$ (consisting of the in-neighbors of $u$ in $\mc G_{P}$) colored in grey.
}
\end{figure}

We conclude this section by the following example showing how Theorem \ref{theorem:main} can be applied even when the perturbed stochastic matrix $\tilde P$ is not irreducible (while the unperturbed one $P$ is). 
\begin{example}\label{example:disconnected}
Let $P$ be an irreducible stochastic matrix with invariant probability distribution $\pi$ and support graph $\mc G_P$ as the one depicted in the leftmost figure above. Consider a node $u\in\mc V$ such that the graph obtained by removing $u$ from $\mc G$ remains strongly connected, and let $\tilde P$ be the stochastic matrix with entries $$\tilde P_{uu}=1\,,\qquad\tilde P_{uv}=\tilde P_{vu}=0\,,\qquad\tilde P_{vv'}=P_{vv'}/(1-P_{vu})\,,$$ for every $v,v'\in\mc V\setminus\{u\}$. Then $\tilde P$ can be interpreted as a perturbation of $P$ with perturbation set $\mc W$ consisting of the node $u$ as well as of all the in-neighbors of $u$ in $\mc G_{P}$, i.e., $\mc W=\{u\}\cup\{v\in\mc V:\,P_{vu}>0\}$. (See Figure \ref{figdisconnected}.) Such setting proves useful when dealing with the removal of a node from, or the addition of a node to, an existing network, e.g., in the context of distributed averaging in sensor networks or opinion dynamics in social networks.

Observe that the matrix $\tilde P$ admits a continuum of invariant probability distributions which is the convex hull of the distribution $\delta^{(u)}$ concentrated in $u$, and of a distribution $\tilde\pi$  supported on $\mc V\setminus\{u\}$. For such $\tilde\pi$, provided that $\mc W\ne\mc V$, we have that $\tilde\gamma_{\mc W}>0$ since the minimization in the righthand side of \eqref{def:exitprobability} runs over all $w\in\mc W\setminus\{u\}$, so that Theorem \ref{theorem:main} can be applied to get a non-trivial upper bound on the total variation distance $||\tilde\pi-\pi||$. On the other hand, for every other invariant probability distribution $\tilde\pi^{(\alpha)}=(1-\alpha)\tilde\pi+\alpha\delta^{(u)}$ for $0<\alpha\le 1$, one gets $\tilde\gamma_{\mc W}=0$, so that Theorem \ref{theorem:main} does not provide any nontrivial bound on $||\tilde\pi^{(\alpha)}-\pi||$.

\end{example}

%
%

\section{Back to the applications}
In this section, we discuss applications of Theorem \ref{theorem:main} first to the PageRank manipulation problem, and then to stochastic matrices associated to networks with a finite-dimensional grid structure. 

\subsection{PageRank manipulation (continued)}\label{pagerank}
For a stochastic matrix $Q$, a probability vector $\mu$, and some $\beta\in(0,1)$, let $P$ and $\pi$ be as in Section \ref{ex:PageRank}.
Let $\tilde Q$ be a perturbation of $Q$, and $\tilde P=(1-\beta)\tilde Q+\beta\1\mu$. Clearly, one has that $\mc W:=\supp(\tilde Q-Q)\supseteq\supp(\tilde P-P)$. 
Moreover, one easily gets the following estimate of the exit probability
\be\label{kappaestimate}\tilde\gamma_{\mc W}\ge\min_{w\in\mc W}\sum_{v\in\mc V\setminus\mc W}P_{wv}\ge\beta(1-\mu(\mc W))\,.\ee
On the other hand, the mixing time can be easily bounded by considering a coupling of two Markov chains, $U(t)$ and $V(t)$ defined as follows. Before meeting, $U(t)$ and $V(t)$ move independently according to the transition probability matrix $Q$ with probability $(1-\beta)$ and jump to a common new state chosen according to $\mu$ with probability $\beta$. Then, starting from the first time they meet, i.e., for $$t\ge T_{c}:=\inf\{t\ge0:\,U(t)=V(t)\}\,,$$ $U(t)=V(t)$ move together with transition probability matrix $P$. 
For every $t\ge0$ and $u,v\in\mc V$,  \cite[Theorem 5.2]{LevinPeresWilmer} implies that
$$||P^t_{u,\cdot}-P^t_{v,\cdot}||_{}\le\P(T_{c}>t|U(0)=u,V(0)=v)\le (1-\beta)^t\,,$$
so that
\be\label{taumixestimate}\tmix\le\l\lceil\frac{-1}{\log(1-\beta)}\r\rceil\le\frac1\beta+1\,.\ee

Finally, let $\tau^{\mu}_{\mc W}:=\sum_{v\in\mc V}\mu_v\tau^v_{\mc W}$ be the expected hitting time of the Markov chain with initial probability distribution $\mu$ and transition probability matrix $P$. 
For all $v$, one has that
$$\tau^v_{\mc W}\le\sum_{t\ge0}(1-\beta)^{t}\beta(t+\tau^{\mu}_{\mc W})=\frac{1-\beta}{\beta}+\tau^{\mu}_{\mc W}\,.$$
Using Kac's formula (\ref{Kac}), the above implies that 
$$\frac1{\pi(\mc W)}=1+\sum_{w\in\mc W}\sum_{v\in\mc V}\frac{\pi_w}{\pi(\mc W)}P_{wv}\tau^v_{\mc W}\le\frac{1}{\beta}+\tau^{\mu}_{\mc W}\,.$$
It follows that 
\be\label{tau*estimate}\tau^*_{\mc W}\ge\beta\tau^{\mu}_{\mc W}\ge\frac{\beta}{\pi(\mc W)}-1 \,.\ee
By combining (\ref{kappaestimate}), (\ref{taumixestimate}), and (\ref{tau*estimate}) with Theorem \ref{theorem:main}, one gets that  
$$||\tilde\pi-\pi||_{}\le\Psi\l(\frac{(1+\beta)\pi(\mc W)}{\beta^2(1-\mu(\mc W))}\r)\,.$$
In particular, the above implies that the alteration of a set of rows $\mc W$ of vanishing aggregate PageRank $\pi(\mc W)$, and $\mu(\mc W)$ bounded away from $1$, has a negligible effect on the whole PageRank vector $\pi$ (in total variation distance).

%
%
%
%
%

\subsection{Networks with high local connectivity}\label{sect:localcoonectivity}
Applications of our results to examples like the distributed averaging algorithm with faulty communication links or to the voter model with influential agents amount to working with perturbations of lazy random walks on graphs, i.e., of stochastic matrices of the form $P=(I+Q)/2$, where $I$ is the identity matrix and $Q$ is the stochastic matrix defined by $Q_{uv}=1/d_u$ if $(u,v)\in\mc E$ and $Q_{uv}=0$ otherwise. The entrance time $\tau^*_{\mc W}$ can be, in general, difficult to be estimated in typical applications when $P$ is sparse and $\mc W$ is a small subset of $\mc V$. 
In this section, we propose some initial results under two assumptions: one is that the set $\mc W$ is not only small but localized in the graph. The second one is that the graph has high local connectivity so that removing $\mc W$ does not drastically alter distances in the remaining part of the graph. The typical graphs for which this holds true are the $d$-dimensional grids (with $d\geq 3$). We believe that both assumptions can be considerably weakened at the price of a deeper analysis. This is the subject of undergoing research which we aim at presenting in another paper.

We start with a simple example to be generalized later on.

\begin{example}\label{ex:torus1}
For integers $m\ge2$ and $d\ge1$, let $P$ be the transition probability matrix of the lazy random walk on a $d$-dimensional toroidal grid of size $n=m^d$. I.e., the node set 
$\mc V=\Z_m^d$ coincides with the direct product of $d$ copies of the cyclic group of integers modulo $m$, and, for all $u,v\in\mc V$, $P_{uu}=1/2$, $P_{uv}=1/(4d)$ if $\sum_{1\le i\le d}|u_i-v_i|=1$, and $P_{uv}=0$ if $\sum_{1\le i\le d}|u_i-v_i|\ge2$. 
For some $w\in\mc V$ and $\alpha\in(0,1)$, consider a perturbed stochastic matrix $\tilde P$ coinciding with $P$ outside $w$, and such that $\tilde P_{ww}<1$. 
Put $\mc W=\{w\}$. It is immediate to verify that 
$$\tilde\gamma_{\mc W}=1-\tilde P_{ww}\,.$$ 
On the other hand, Kac's formula (\ref{Kac}) implies that 
$$n=\frac1{\pi_w}=1+\frac1{4d}\sum_{v:|v-w|=1}\tau^v_w=1+\frac12\tau^*_{\mc W}\,,$$ 
where last equality follows from a basic symmetry argument. 
Moreover, standard results \cite[Theorem 5.5]{LevinPeresWilmer} imply that 
$$\tmix\le C_d n^{2/d}$$
for some constant $C_d$ depending on $d$ but not on $n$. 
Then, Theorem \ref{theorem:main} implies that 
$$||\pi-\tilde\pi||_{}\le\Psi\l(\frac{2C_d}{1-\tilde P_{ww}}\,\cdot\,\frac{n^{2/d}}{n-1} \r)\,.$$
The above guarantees that the total variation distance $||\pi-\tilde\pi||_{}$ vanishes as $n$ grows large provided that $d\ge3$.  
\end{example}


In the previous example, $\tau^*_{\mc W}$ was exactly computed in terms of $\pi(\mc W)$ using Kac's formula \eqref{Kac} and the spatial symmetry in the neighborhood of the perturbed set $\mc W=\{w\}$. For general $\mc W$, such symmetry argument breaks down. Below we propose a way to overcome this difficulty 
in a general situation where  $\mc W$ is localized and its boundary is sufficiently well connected in $\mc V\setminus\mc W$.
\begin{figure}\begin{center}
{\includegraphics[height=4.5cm,width=6.5cm]{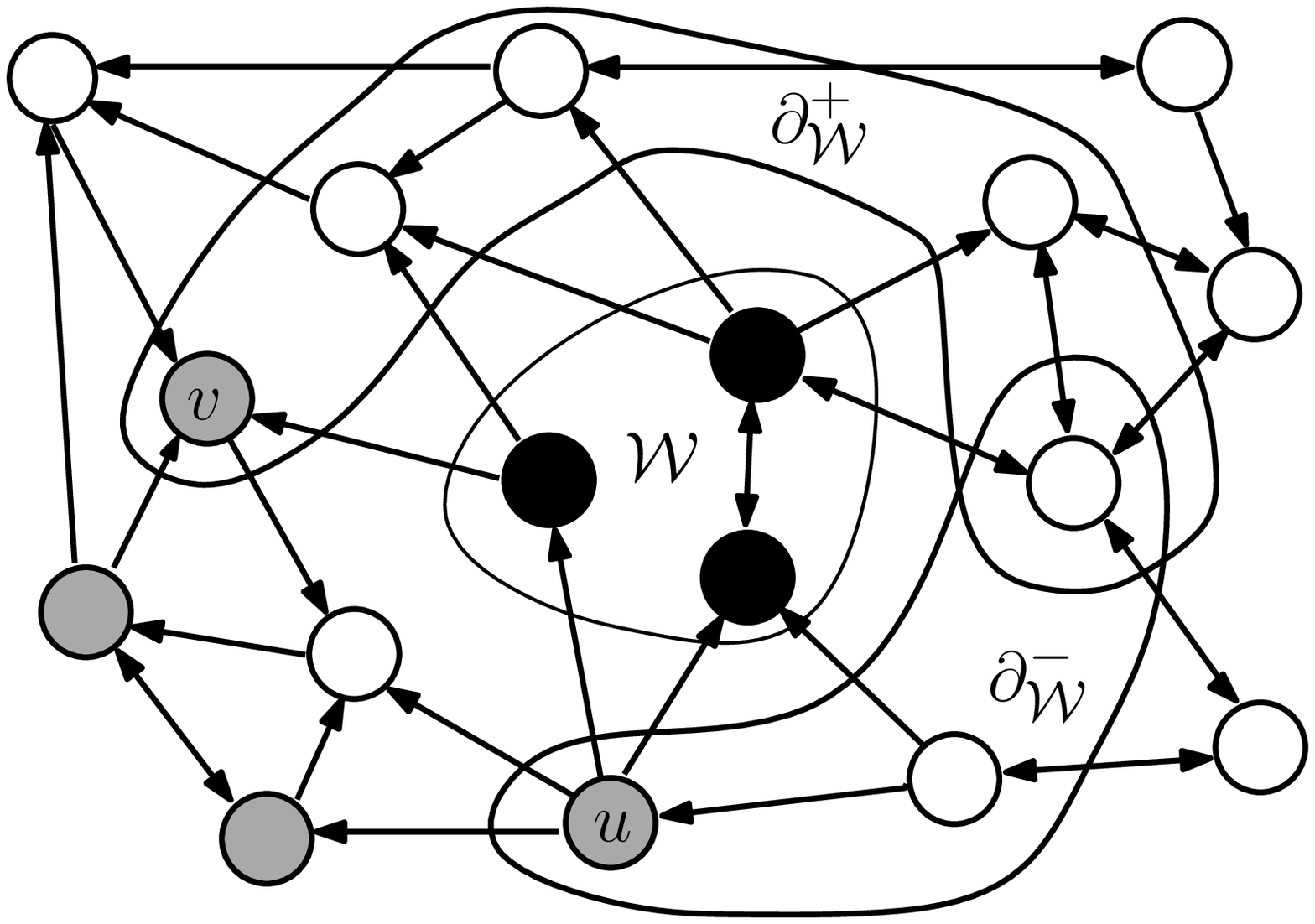}}
\caption{\label{fig:WdeW} The external boundaries $\partial^+_{\mc W}$ and $\partial^-_{\mc W}$ of a node set $\mc W$. A simple path in $\mc V\setminus\mc W$ from $u\in\partial^-_{\mc W}$ to $v\in\partial^+_{\mc W}$ is shaded in gray.}
\end{center}
\end{figure}
Define the external boundaries of $\mc W$ as
$$\partial^+_{\mc W}:=\{v\in \mc V\setminus\mc W\,:\  P_{wv}>0\text{ for some }w\in\mc W\}\,,$$
$$\partial^-_{\mc W}:=\{v\in \mc V\setminus\mc W\,:\  P_{vw}>0\text{ for some }w\in\mc W\}\,.$$
(See Figure \ref{fig:WdeW}.) Clearly,
\be\label{tau*property}\tau^*_{\mc W}=\min\{\tau^v_{\mc W}:\,v\in\partial^-_{\mc W}\}\,.\ee
On the the other hand, let 
\be\label{taucircdef}\tau^{\circ}_{\mc W}:=\max\{\tau^v_{\mc W}\,:\ v\in\partial^+_{\mc W}\}\,,\ee
and observe that, from Kac's formula (\ref{Kac}),
\begin{equation}\label{maxhit} \tau^{\circ}_{\mc W}\geq \sum_{w\in\mc W}\sum_{v\in\mc V}\frac{\pi_w}{\pi(\mc W)} P_{wv}\tau^v_{\mc W}=\frac1{\pi(\mc W)}-1\,.\end{equation}
Now, for all $u\in\partial^-_{\mc W}$ and $v\in\partial^+_{\mc W}$,  let $\Gamma_{u,v}$ be the (possibly empty) set of simple paths in $ \mc V\setminus\mc W$ starting in $u$ and ending in $v$. 
For all paths $\xi=(u=\xi_0,\xi_1,\ldots,\xi_l=v)\in\Gamma_{u,v}$,  let $P_{\xi}:=\prod_{1\le i\le l}P_{\xi_{i-1}\xi_i}$. Define 
\begin{equation}\label{connect}\lambda_{\mc W}:=\min_{u,v}\max_{\xi\in\Gamma_{u,v}}P_{\xi}\,,\ee
where the minimization is intended to run over all $u\in\partial^+_{\mc W}$ and $v\in\partial^-_{\mc W}$ such that $u\ne v$, and we use the convention that the minimum over an empty set equals $1$, and the maximum over an empty set equals $0$. 
Then, the following result holds true. 
\begin{proposition}\label{lemma: well-connected}
Let $P$ be an irreducible stochastic matrix on a finite set $\mc V$, and $\pi=\pi P$ its invariant probability vector. Then, for all $\mc W\subseteq\mc V$, the entrance time $\tau^*_{\mc W}$ satisfies
$$\tau^*_{\mc W}\geq  \lambda_{\mc W}\l(\frac1{\pi(\mc W)}-1\right)\,,$$
where $\lambda_{\mc W}$ is defined as in \eqref{connect}.
\end{proposition}
\proof
Let $u\in\partial^-_{\mc W}$ and $v\in\partial^+_{\mc W}$ be such that $\tau^{u}_{\mc W}=\tau^{*}_{\mc W}$ and $\tau^{v}_{\mc W}=\tau^{\circ}_{\mc W}$.
For a path $\xi=(u=\xi_0,\xi_1,\ldots,\xi_{l-1},\xi_l=v)$ in $\Gamma_{u,v}$, let $\1_{\xi}$ be the indicator function of the event that the first $l$ steps of the Markov chain $V(t)$ started at $u$ and moving with transition probability matrix $P$ are along $\xi$, i.e., the event $\cap_{t=0}^l\{V(t)=\xi_t\}$. Then,
\begin{equation}\label{minhit}
\tau^*_{\mc W}=\tau^u_{\mc W}\geq \E_u[T_{\mc W}\1_{\xi}]=P_{\xi}(\tau^{v}_{\mc W}+l)\ge P_{\xi}\tau^{\circ}_{\mc W}\,.\end{equation}
The claim now follows from (\ref{maxhit}), (\ref{minhit}), and the arbitrariness of $\xi\in\Gamma_{u,v}$.
\qed

The above result turns out to be useful in those contexts where the set $\mc W$ is sufficiently localized so that its boundary is tightly connected outside $\mc W$ and $\lambda_{\mc W}$ remains bounded away from $0$. 

\begin{example}\label{ex:torus3}
Let $P$ be the lazy simple random walk on a $d$-dimensional toroidal grid as in Example \ref{ex:torus1} 
and let the perturbation be supported on a  hypercube $\mc W=\prod_{i=1}^d[\alpha_i,\alpha_i+s-1]$. 
One can easily verify that any pair of nodes in $\partial^+_{\mc W}=\partial^-_{\mc W}$ can be connected by a path of length $d(s+1)$ outside $\mc W$, so that $\lambda_{\mc W}\geq (4d)^{-d(s+1)}$. On the other hand, $n\pi(\mc W)=|\mc W|=s^d$, so that Proposition \ref{lemma: well-connected} implies that 
 $$\tau^*_{\mc W}\geq \frac1{(4d)^{d(s+1)}}\left(\frac1{\pi(\mc W)}-1\right)= \frac1{(4d)^{d(s+1)}}\left(\frac n{s^{d}}-1\right)\,.$$ Since the mixing time satisfies $$\tmix \leq C_dn^{2/d}$$ for some positive constant $C_d$ independent from $n$ \cite[Theorem 5.5]{LevinPeresWilmer}, we have that 
\be\label{est1}\frac{\tmix}{\tau^*_{\mc W}}\leq C'_d\frac{(4d)^{ds}  s^d}{1-  s^d/n}n^{2/d-1}\,,\ee
with $C'_d:=C_d(4d)^d$.

It remains to be estimated the exit probability from
 $\mc W$ which is the (only) term depending on finer details of the perturbed matrix $\tilde P$. Assume that $\tilde P$ is irreducible, and put 
$$\delta =\min\l\{\tilde P_{wv}\,:\, w\in \mc W\,,\; \tilde P_{wv}>0\r\}\,.$$
Since from every $w\in \mc W$ there is a path leading to $\partial\mc W$ of length at most $|\mc W|=s^d$, one gets that 
\be\label{est2}\tilde\gamma_{\mc W}\geq \min_{\substack{w:\,\tilde\pi_w>0}}\frac{1}{s^d}\P_w(\tilde T_{\mc V\setminus\mc W}\leq s^d)\geq \frac{\delta^{s^d}}{s^d}\,.\ee
By dividing both sides of \eqref{est1} by the respective sides of \eqref{est2}, one gets 
$$
\frac{\tmix}{\tilde\gamma_{\mc W}\cdot\tau^*_{\mc W}}\leq C'_d\frac{(4d)^{ds}s^{2d}\delta^{-s^d}}{(1-  s^d/n)n^{1-2/d}}\,.
$$
Note that the righthand side of the above vanished as $n$ grows large provided that 
the term $s^{2d}\delta^{-s^d}$ grows slower than $n^{1-2/d}$.  
Then, Theorem \ref{theorem:main} implies  that a sufficient condition for  
$||\tilde\pi-\pi||\to 0$ as $n$ grows large is that $|\mc W|=s^d$ satisfies 
$$\limsup_n\frac{|\mc W|}{\log n}<\frac{d-2}{d\log \delta^{-1}}\,.$$

\end{example}

\section{Conclusion}
Invariant probability vectors of stochastic matrices play a central role in a large number of multi-agent network problems including distributed averaging algorithms, opinion dynamics, and centrality measures such as PageRank. This paper investigates the fundamental question of how resilient such invariant probability vectors are with respect to perturbations of the network. The main result provides an estimate of the total variation distance between the invariant probability vectors of two stochastic matrices in terms of the mixing time of one of the matrices and of the size of the perturbation set $\mc W$ measured as the product of two quantities: the entrance time on $\mc W$ and the exit probability from $\mc W$. Explicit applications to network models have also been discussed in detail. Among the relevant issues which have not been addressed by this paper and deserve to be considered for future research are:

\begin{itemize}
\item The estimation of the entrance time of the perturbation set remains the most challenging problem in applying our result. In particular, we would like to extend our estimation to small but scattered perturbation sets as well to other general classes of networks such as locally tree-like graphs.

\item In many applications of network centrality, the total variation distance between two probability vectors may not be the most relevant measure of the effect of a perturbation. E.g., the maximal ratio of the centralities assigned to the same node in the unperturbed and in the perturbed network would be of great potential interest in such cases.

\item When a network is perturbed locally, we expect the effect of the perturbation to decay as a function of the distance from the perturbation set. This is not captured by the total variation analysis and may require an essentially different approach.

\end{itemize}

\section*{Acknowledgements}
The first author is a member of the LCCC Linnaeus Center and the ELLIIT Excellence Center at Lund University. His research has been partially supported by the Swedish Research Council (VR) through the the junior research grant Information Dynamics in Large-Scale Networks.

\bibliographystyle{amsplain}
\bibliography{influence}

\begin{IEEEbiography}[{\includegraphics[width=1in,height=1.25in,clip,keepaspectratio]{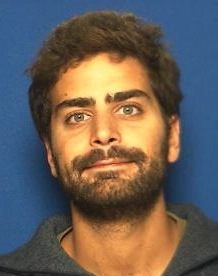}}]{Giacomo Como}
 (M'12) is an Associate Professor at the Department of Automatic Control, Lund University, Sweden. He received the B.Sc., M.S., and Ph.D. degrees in Applied Mathematics from Politecnico di Torino, Italy, in 2002, 2004, and 2008, respectively. In 2006--2007, he was a Visiting Assistant in Research at the Department of Electrical Engineering, Yale University. From 2008 to 2011, he was a Postdoctoral Associate at the Laboratory for Information and Decision
Systems, Massachusetts Institute of Technology. His current research interests are in control, information, game theory, and networks.
\end{IEEEbiography}

\begin{IEEEbiography}[{\includegraphics[width=1in,height=1.25in,clip,keepaspectratio]{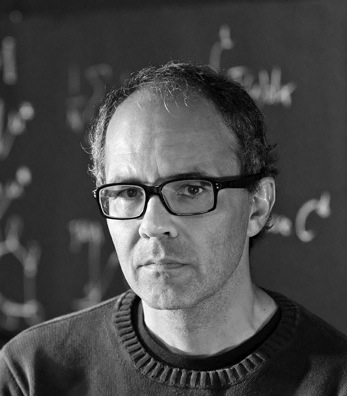}}]{Fabio Fagnani }
received the Laurea degree in mathematics from the University of Pisa, Pisa, Italy, and from Scuola Normale Superiore of Pisa, Pisa, Italy,
in 1986. He got the Ph.D.~degree in mathematics from the University of Groningen, Groningen, the Netherlands, in 1991.
From 1991 to 1998, he was an Assistant Professor of mathematical analysis with Scuola Normale Superiore.
In 1997, he was a Visiting Professor with Massachusetts Institute of Technology (MIT), Cambridge, MA, USA. Since 1998, he has been with the Politecnico of Torino, where he is currently (since 2002) a Full Professor of mathematical analysis. From 2006 to 2012, he has acted as a Coordinator of the Ph.D.~program Mathematics for Engineering Sciences of Politecnico di Torino, Torino, Italy. From June 2012, he is the Head of the Department of Mathematical Sciences, Politecnico di Torino. His current research topics are on cooperative algorithms and dynamical systems over graphs, inferential distributed algorithms, and opinion dynamics. He has published over 50 refereed papers on international journals, he has delivered invited conferences in many international workshops and conferences and in many universities. 
He is an Associate Editor of the IEEE Transactions on Network Science and Engineering and of the IEEE Transactions on Control of Network Systems. He is a Member of the International Program Committee for the events NECSYS.
\end{IEEEbiography}
\end{document}